\documentclass[times,10pt,twocolumn]{article} 
\usepackage{latex8}
\usepackage{times}
\usepackage{amsmath}
\usepackage{amsthm}
\usepackage{amsfonts}
\usepackage{amssymb}

\newcommand\R{{\mathbf{R}}}
\newcommand\Z{{\mathbf{Z}}}
\newcommand\F{{\mathbf{F}}}
\newcommand\E{{\mathbf{E}}}
\newcommand\C{{\mathbf{C}}}
\newcommand\X{{\mathbf{X}}}
\newcommand\Y{{\mathbf{Y}}}
\newcommand\B{{\mathcal{B}}}
\renewcommand\S{{\mathcal{S}}}
\newcommand\eps{{\varepsilon}}
\newcommand\str{{\operatorname{str}}}
\newcommand\psd{{\operatorname{psd}}}
\newcommand\err{{\operatorname{err}}}

\theoremstyle{plain}
  \newtheorem{theorem}[subsection]{Theorem}
  \newtheorem{conjecture}[subsection]{Conjecture}

  \newtheorem{proposition}[subsection]{Proposition}
  
  \newtheorem{lemma}[subsection]{Lemma}
  \newtheorem{corollary}[subsection]{Corollary}
   
\theoremstyle{remark}
  \newtheorem{remark}[subsection]{Remark}
  
   \newtheorem{example}[subsection]{Example}
\theoremstyle{definition}
  \newtheorem{definition}[subsection]{Definition}

\begin{document}

\title{Structure and randomness in combinatorics}

\author{Terence Tao\\
Department of Mathematics, UCLA\\405 Hilgard Ave, Los Angeles CA 90095\\tao@math.ucla.edu
}

\maketitle
\thispagestyle{empty}

\begin{abstract}
   Combinatorics, like computer science, often has to deal with large objects of unspecified (or unusable) structure.   One powerful way to deal with such an arbitrary object is to decompose it into more usable components.  In particular, it has proven profitable to decompose such objects into a \emph{structured} component, a \emph{pseudo-random} component, and a \emph{small} component (i.e. an error term); in many cases it is the structured component which then dominates.  We illustrate this philosophy in a number of model cases.
\end{abstract}

\Section{Introduction}

In many situations in combinatorics, one has to deal with an object of large complexity or entropy - such as a graph on $N$ vertices, a function on $N$ points, etc., with $N$ large.  We are often interested in the worst-case behaviour of such objects; equivalently, we are interested in obtaining results which apply to \emph{all} objects in a certain class, as opposed to results for almost all objects (in particular, random or average case behaviour) or for very specially structured objects.  The difficulty here is that the spectrum of behaviour of an arbitrary large object can be very broad.  At one extreme, one has very \emph{structured} objects, such as complete bipartite graphs, or functions with periodicity, linear or polynomial phases, or other algebraic structure.  At the other extreme are \emph{pseudorandom} objects, which mimic the behaviour of random objects in certain key statistics (e.g. their correlations with other objects, or with themselves, may be close to those expected of random objects).

Fortunately, there is a fundamental phenomenon that one often has a \emph{dichotomy} between structure and pseudorandomness, in that given a reasonable notion of structure (or pseudorandomness), there often exists a dual notion of pseudorandomness (or structure) such that an arbitrary object can be decomposed into a structured component and a pseudorandom component (possibly with a small error).  Here are two simple examples of such decompositions:

\begin{itemize}
\item[(i)] An \emph{orthogonal decomposition} $f = f_\str + f_\psd$ of a vector $f$ in a Hilbert space into its orthogonal projection $f_\str$ onto a subspace $V$ (which represents the ``structured'' objects), plus its orthogonal projection $f_{\psd}$ onto the orthogonal complement $V^\perp$ of $V$ (which represents the ``pseudorandom'' objects).  
\item[(ii)] A \emph{thresholding} $f = f_\str + f_\psd$ of a vector $f$, where $f$ is expressed in terms of some basis $v_1,\ldots,v_n$ (e.g. a Fourier basis) as $f = \sum_{1 \leq i \leq n} c_i v_i$, the ``structured'' component $f_\str := \sum_{i: |c_i| \geq \lambda} c_i v_i$ contains the contribution of the large coefficients, and the ``pseudorandom'' component $f_\psd := \sum_{i: |c_i| < \lambda} c_i v_i$ contains the contribution of the small coefficients.  Here $\lambda > 0$ is a thresholding parameter which one is at liberty to choose.
\end{itemize}

Indeed, many of the decompositions we discuss here can be viewed as variants or perturbations of these two simple decompositions.  More advanced examples of decompositions include the Szemer\'edi regularity lemma for graphs (and hypergraphs), as well as various \emph{structure theorems} relating to the Gowers uniformity norms, used for instance in \cite{gt-primes}, \cite{gt-finite}.  Some decompositions from classical analysis, most notably the \emph{spectral decomposition} of a self-adjoint operator into orthogonal subspaces associated with the pure point, singular continuous, and absolutely continuous spectrum, also have a similar spirit to the structure-randomness dichtomy.

The advantage of utilising such a decomposition is that one can use different techniques to handle the structured component and the pseudorandom component (as well as the error component, if it is present).  Broadly speaking, the structured component is often handled by algebraic or geometric tools, or by reduction to a ``lower complexity'' problem than the original problem, whilst the contribution of the pseudorandom and error components is shown to be negligible by using inequalities from analysis (which can range from the humble Cauchy-Schwarz inequality to other, much more advanced, inequalities).  A particularly notable use of this type of decomposition occurs in the many different proofs of Szemer\'edi's theorem \cite{szemeredi}; see e.g. \cite{tao:montreal} for further discussion. 

In order to make the above general strategy more concrete, one of course needs to specify more precisely what ``structure'' and ``pseudorandomness'' means.  There is no single such definition of these concepts, of course; it depends on the application.  In some cases, it is obvious what the definition of one of these concepts is, but then one has to do a non-trivial amount of work to describe the dual concept in some useful manner.  We remark that \emph{computational} notions of structure and randomness do seem to fall into this framework, but thus far all the applications of this dichotomy have focused on much simpler notions of structure and pseudorandomness, such as those associated to Reed-Muller codes.

In these notes we give some illustrative examples of this structure-randomness dichotomy.  While these examples are somewhat abstract and general in nature, they should by no means be viewed as the definitive expressions of this dichotomy; in many applications one needs to modify the basic arguments given here in a number of ways.  On the other hand, the core ideas in these arguments (such as a reliance on energy-increment or energy-decrement methods) appear to be fairly universal.  The emphasis here will be on illustrating the ``nuts-and-bolts'' of structure theorems; we leave the discussion of the more advanced structure theorems and their applications to other papers.

One major topic we will not be discussing here (though it is lurking underneath the surface) is the role of ergodic theory in all of these decompositions; we refer the reader to \cite{tao:montreal} for further discussion.  Similarly, the recent ergodic-theoretic approaches to hypergraph regularity, removal, and property testing in \cite{tao:infinite}, \cite{austin} will not be discussed here, in order to prevent the exposition from becoming too unfocused.  The lecture notes here also have some intersection with the author's earlier article \cite{tao-icm}.

\Section{Structure and randomness in a Hilbert space}

Let us first begin with a simple case, in which the objects one is studying lies in some real finite-dimensional Hilbert space $H$, and the concept of structure is captured by some known set $S$ of ``basic structured objects''.  This setting is already strong enough to establish the Szemer\'edi regularity lemma, as well as variants such as Green's arithmetic regularity lemma.  One should think of the dimension of $H$ as being extremely large; in particular, we do not want any of our quantitative estimates to depend on this dimension.

More precisely, let us designate a finite collection $S \subset H$ of ``basic structured'' vectors of bounded length; we assume for concreteness that $\|v\|_H \leq 1$ for all $v \in S$.  We would like to view elements of $H$ which can be ``efficiently represented'' as linear combinations of vectors in $S$ as \emph{structured}, and vectors which have low correlation (or more precisely, small inner product) to all vectors in $S$ as \emph{pseudorandom}.  More precisely, given $f \in H$, we say that $f$ is \emph{$(M,K)$-structured} for some $M, K > 0$ if one has a decomposition
$$ f = \sum_{1 \leq i \leq M} c_i v_i$$
with $v_i \in S$ and $c_i \in [-K,K]$ for all $1 \leq i \leq M$.  We also say that $f$ is \emph{$\eps$-pseudorandom} for some $\eps > 0$ if we have $|\langle f,v \rangle_H| \leq \eps$ for all $v \in S$.  It is helpful to keep some model examples in mind:

\begin{example}[Fourier structure]\label{fourier-struct}  Let $\F_2^n$ be a Hamming cube; we identify the finite field $\F_2$ with $\{0,1\}$ in the usual manner.  We let $H$ be the $2^n$-dimensional space of functions $f: \F_2^n \to \R$, endowed with the inner product
$$ \langle f,g \rangle_H := \frac{1}{2^n} \sum_{x \in \F_2^n} f(x) g(x),$$
and let $S$ be the space of \emph{characters},
$$ S := \{ e_\xi: \xi \in \F_2^n\},$$
where for each $\xi \in \F_2^n$, $e_\xi$ is the function $e_\xi(x) := (-1)^{x \cdot \xi}$.  Informally, a structured function $f$ is then one which can be expressed in terms of a small number (e.g. $O(1)$) characters, whereas a pseudorandom function $f$ would be one whose Fourier coefficients \begin{equation}\label{fourier-def}
\hat f(\xi) := \langle f, e_\xi \rangle_H
\end{equation}
 are all small.
\end{example}

\begin{example}[Reed-Muller structure]  Let $H$ be as in the previous example, and let $1 \leq k \leq n$.  We now let $S = S_k(\F_2^n)$ be the space of Reed-Muller codes $(-1)^{P(x)}$, where $P: \F_2^n \to \F_2$ is any polynomial of $n$ variables with coefficients and degree at most $k$.  For $k=1$, this gives the same notions of structure and pseudorandomness as the previous example, but as we increase $k$, we enlarge the class of structured functions and shrink the class of pseudorandom functions.  For instance, the function $(x_1,\ldots,x_n) \mapsto (-1)^{\sum_{1 \leq i < j \leq n} x_i x_j}$ would be considered highly pseudorandom when $k=1$ but highly structured for $k \geq 2$.
\end{example}

\begin{example}[Product structure]\label{prod}  Let $V$ be a set of $|V|=n$ vertices, and let $H$ be the $n^2$-dimensional space of functions $f: V \times V \to \R$, endowed with the inner product
$$ \langle f, g \rangle_H := \frac{1}{n^2} \sum_{v,w \in V} f(v,w) g(v,w).$$
Note that any graph $G = (V,E)$ can be identified with an element of $H$, namely the indicator function $1_E: V \times V \to \{0,1\}$ of the set of edges.  We let $S$ be the collection of tensor products $(v,w) \mapsto 1_A(v) 1_B(w)$, where $A, B$ are subsets of $V$.  Observe that $1_E$ will be quite structured if $G$ is a complete bipartite graph, or the union of a bounded number of such graphs.  At the other extreme, if $G$ is an \emph{$\eps$-regular} graph of some edge density $0 < \delta < 1$ for some $0 < \eps < 1$, in the sense that the number of edges between $A$ and $B$ differs from $\delta |A| |B|$ by at most $\eps |A| |B|$ whenever $A, B \subset V$ with $|A|, |B| \geq \eps n$, then $1_E - \delta$ will be $O(\eps)$-pseudorandom.
\end{example}

We are interested in obtaining quantative answers to the following general problem: given an arbitrary bounded element $f$ of the Hilbert space $H$ (let us say $\|f\|_H \leq 1$ for concreteness), can we obtain a decomposition
\begin{equation}\label{fdecomp}
f = f_\str + f_\psd + f_\err
\end{equation}
where $f_{\str}$ is a structured vector, $f_\psd$ is a pseudorandom vector, and $f_\err$ is some small error?

One obvious ``qualitative'' decomposition arises from using the vector space $\operatorname{span}(S)$ spanned by the basic structured vectors $S$.  If we let $f_\str$ be the orthogonal projection from $f$ to this vector space, and set $f_\psd := f - f_\str$ and $f_\err := 0$, then we have perfect control on the pseudorandom and error components: $f_\psd$ is $0$-pseudorandom and $f_\err$ has norm $0$.  On the other hand, the only control on $f_\str$ we have is the qualitative bound that it is $(K,M)$-structured for some finite $K, M < \infty$.  In the three examples given above, the vectors $S$ in fact span all of $H$, and this decomposition is in fact trivial!

We would thus like to perform a tradeoff, increasing our control of the structured component at the expense of worsening our control on the pseudorandom and error components.  We can see how to achieve this by recalling how the orthogonal projection of $f$ to $\operatorname{span}(S)$ is actually constructed; it is the vector $v$ in $\operatorname{span}(S)$ which minimises the ``energy'' $\|f-v\|_H^2$ of the residual $f-v$.  The key point is that if $v \in \operatorname{span}(S)$ is such that $f-v$ has a non-zero inner product with a vector $w \in S$, then it is possible to move $v$ in the direction $w$ to decrease the energy $\|f-v\|_H^2$.  We can make this latter point more quantitative:

\begin{lemma}[Lack of pseudorandomness implies energy decrement]\label{pseud-dec}  Let $H, S$ be as above.  Let $f \in H$ be a vector with $\|f\|_H^2 \leq 1$, such that $f$ is \emph{not} $\eps$-pseudorandom for some $0 < \eps \leq 1$.  Then there exists $v \in S$ and $c \in [-1/\eps,1/\eps]$ such that $|\langle f, v \rangle| \geq \eps$ and
$\|f - cv \|_H^2 \leq \|f\|_H^2 - \eps^2$.
\end{lemma}

\begin{proof}  By hypothesis, we can find $v \in S$ be such that $|\langle f, v \rangle| \geq \eps$, thus by Cauchy-Schwarz and hypothesis on $S$
$$ 1 \geq \|v\|_H \geq |\langle f, v \rangle| \geq \eps.$$  
We then set $c := \langle f, v \rangle / \|v\|_H^2$ (i.e. $cv$ is the orthogonal projection of $f$ to the span of $v$).  The claim then follows from Pythagoras' theorem.
\end{proof}

If we iterate this by a straightforward greedy algorithm argument we now obtain

\begin{corollary}[Non-orthogonal weak structure theorem]\label{wst}  Let $H, S$ be as above.  Let $f \in H$ be such that $\|f\|_H \leq 1$, and let $0 < \eps \leq 1$.  Then there exists a decomposition \eqref{fdecomp} such that $f_\str$ is $(1/\eps^2, 1/\eps)$-structured, $f_\psd$ is $\eps$-pseudorandom, and $f_\err$ is zero.
\end{corollary}

\begin{proof}  We perform the following algorithm.
\begin{itemize}
\item Step 0.  Initialise $f_\str := 0$, $f_\err := 0$, and $f_\psd := f$.  Observe that $\|f_\psd\|_H^2 \leq 1$.
\item Step 1.  If $f_\psd$ is $\eps$-pseudorandom then {\textbf STOP}.  Otherwise, by Lemma \ref{pseud-dec}, we can find $v \in S$ and $c \in [-1/\eps,1/\eps]$ such that $\|f_\psd - cv \|_H^2 \leq \|f_\psd\|_H^2 - \eps^2$.
\item Step 2.  Replace $f_\psd$ by $f_\psd - cv$ and replace $f_\str$ by $f_\str + cv$.  Now return to Step 1.
\end{itemize}

It is clear that the ``energy'' $\|f_\psd\|_H^2$ decreases by at least $\eps^2$ with each iteration of this algorithm, and thus this algorithm terminates after at most $1/\eps^2$ such iterations.  The claim then follows.
\end{proof}

Corollary \ref{wst} is not very useful in applications, because the control on the structure of $f_\str$ are relatively poor compared to the pseudorandomness of $f_\psd$ (or vice versa).  One can do substantially better here, by allowing the error term $f_\err$ to be non-zero.  More precisely, we have

\begin{theorem}[Strong structure theorem]\label{sst}  Let $H, S$ be as above, let $\eps > 0$, and let $F: \Z^+ \to \R^+$ be an arbitrary function.  Let $f \in H$ be such that $\|f\|_H \leq 1$. Then we can find an integer $M = O_{F,\eps}(1)$ and a decomposition \eqref{fdecomp} where $f_\str$ is $(M,M)$-structured, $f_\psd$ is $1/F(M)$-pseudorandom, and $f_\err$ has norm at most $\eps$.  
\end{theorem}

Here and in the sequel, we use subscripts in the $O()$ asymptotic notation to denote that the implied constant depends on the subscripts.  For instance, $O_{F,\eps}(1)$ denotes a quantity bounded by $C_{F,\eps}$, for some quantity $C_{F,\eps}$ depending only on $F$ and $\eps$.  Note that the pseudorandomness of $f_\psd$ can be of arbitrarily high quality compared to the complexity of $f_\str$, since we can choose $F$ to be whatever we please; the cost of doing so, of course, is that the upper bound on $M$ becomes worse when $F$ is more rapidly growing.

To prove Theorem \ref{sst}, we first need a variant of Corollary \ref{wst} which gives some orthogonality between $f_\str$ and $f_\psd$, at the cost of worsening the complexity bound on $f_\str$.

\begin{lemma}[Orthogonal weak structure theorem]\label{wst2}  Let $H, S$ be as above.  Let $f \in H$ be such that $\|f\|_H \leq 1$, and let $0 < \eps \leq 1$.  Then there exists a decomposition \eqref{fdecomp} such that $f_\str$ is $(1/\eps^2, O_\eps(1))$-structured, $f_\psd$ is $\eps$-pseudorandom, $f_\err$ is zero, and $\langle f_\str, f_\psd \rangle_H = 0$.
\end{lemma}

\begin{proof}  We perform a slightly different iteration to that in Corollary \ref{wst}, where we insert an additional orthogonalisation step within the iteration to a subspace $V$:
\begin{itemize}
\item Step 0.  Initialise $V := \{0\}$ and $f_\err := 0$.
\item Step 1.  Set $f_\str$ to be the orthogonal projection of $f$ to $V$, and $f_\psd := f - f_\str$. \item Step 2.  If $f_\psd$ is $\eps$-pseudorandom then {\textbf STOP}.  Otherwise, by Lemma \ref{pseud-dec}, we can find $v \in S$ and $c \in [-1/\eps,1/\eps]$ such that $|\langle f_\psd, v \rangle_H| \geq \eps$ and $\|f_\psd - cv \|_H^2 \leq \|f_\psd\|_H^2 - \eps^2$.
\item Step 3.  Replace $V$ by $\operatorname{span}(V \cup \{v\})$, and return to Step 1.
\end{itemize}

Note that at each stage, $\|f_\psd\|_H$ is the minimum distance from $f$ to $V$.  Because of this, we see that $\|f_\psd\|_H^2$ decreases by at least $\eps^2$ with each iteration, and so this algorithm terminates in at most $1/\eps^2$ steps.

Suppose the algorithm terminates in $M$ steps for some $M \leq 1/\eps^2$.  Then we have constructed a nested flag
$$ \{0\} = V_0 \subset V_1 \subset \ldots \subset V_M$$
of subspaces, where each $V_i$ is formed from $V_{i-1}$ by adjoining a vector $v_i$ in $S$.  Furthermore, by construction we have $|\langle f_i, v_i \rangle| \geq \eps$ for some vector $f_i$ of norm at most $1$ which is orthogonal to $V_{i-1}$.  Because of this, we see that $v_i$ makes an angle of $\Theta_\eps(1)$ with $V_{i-1}$.  As a consequence of this and the Gram-Schmidt orthogonalisation process, we see that $v_1,\ldots,v_i$ is a well-conditioned basis of $V_i$, in the sense that any vector $w \in W_i$ can be expressed as a linear combination of $v_1,\ldots,v_i$ with coefficients of size $O_{\eps,i}(\|w\|_H)$.  In particular, since $f_\str$ has norm at most $1$ (by Pythagoras' theorem) and lies in $V_M$, we see that $f_\str$ is a linear combination of $v_1,\ldots,v_M$ with coefficients of size $O_{M,\eps}(1) = O_\eps(1)$, and the claim follows.
\end{proof}

We can now iterate the above lemma and use a pigeonholing argument to obtain the strong structure theorem.

\begin{proof}[Proof of Theorem \ref{sst}]
We first observe that it suffices to prove a weakened version of Theorem \ref{sst} in which $f_\str$ is $(O_{M,\eps}(1), O_{M,\eps}(1))$-structured rather than $(M,M)$ structured.  This is because one can then recover the original version of Theorem \ref{sst} by making $F$ more rapidly growing, and redefining $M$; we leave the details to the reader.  Also, by increasing $F$ if necessary we may assume that $F$ is integer-valued and $F(M) > M$ for all $M$.

We now recursively define $M_0 := 1$ and $M_i := F(M_{i-1})$ for all $i \geq 1$.  We then recursively define $f_0, f_1, \ldots$ by setting $f_0 := f$, and then for each $i \geq 1$ using Lemma \ref{wst2} to decompose $f_{i-1} = f_{\str,i} + f_i$ where $f_{\str,i}$ is $(O_{M_i}(1),O_{M_i}(1))$-structured, and $f_i$ is $1/M_i$-pseudorandom and orthogonal to $f_{\str,i}$.  From Pythagoras' theorem we see that the quantity $\|f_i\|_H^2$ is decreasing, and varies between $0$ and $1$.  By the pigeonhole principle, we can thus find $1 \leq i \leq 1/\eps^2+1$ such that $\|f_{i-1}\|_H^2 - \|f_i\|_H^2 \leq \eps^2$; by Pythagoras' theorem, this implies that $\|f_{\str,i}\|_H \leq \eps$.  If we then set $f_\str := f_{\str,0} + \ldots + f_{\str,i-1}$, $f_\psd := f_i$, $f_\err := f_{\str,i}$, and $M := M_{i-1}$, we obtain the claim.
\end{proof}

\begin{remark} By tweaking the above argument a little bit, one can also ensure that the quantities $f_\str, f_\psd, f_\err$ in Theorem \ref{sst} are orthogonal to each other.  We leave the details to the reader.
\end{remark}

\begin{remark} The bound $O_{F,\eps}(1)$ on $M$ in Theorem \ref{sst} is quite poor in practice; roughly speaking, it is obtained by iterating $F$ about $O(1/\eps^2)$ times.  Thus for instance if $F$ is of exponential growth (which is typical in applications), $M$ can be tower-exponential size in $\eps$.  These excessively large values of $M$ unfortunately seem to be necessary in many cases, see e.g. \cite{gowers-sz} for a discussion in the case of the Szemer\'edi regularity lemma, which can be deduced as a consequence of Theorem \ref{sst}.
\end{remark}

To illustrate how the strong regularity lemma works in practice, we use it to deduce the arithmetic regularity lemma of Green \cite{green} (applied in the model case of the Hamming cube $\F_2^n$).  Let $A$ be a subset of $\F_2^n$, and let $1_A: \F_2^n \to \{0,1\}$ be the indicator function.  If $V$ is an affine subspace (over $\F_2$) of $\F_2^n$, we say that $A$ is \emph{$\eps$-regular} in $V$ for some $0 < \eps < 1$ if we have
$$ \left|\E_{x \in V} (1_A(x) - \delta_V) e_\xi(x)\right| \leq \eps$$
for all characters $e_\xi$, where $\E_{x \in V} f(x) := \frac{1}{|V|} \sum_{x \in V} f(x)$ denotes the average value of $f$ on $V$, and $\delta_V := \E_{x \in V} 1_A(x) = |A \cap V|/|V|$ denotes the density of $A$ in $V$.  The following result is analogous to the celebrated Szemer\'edi regularity lemma:

\begin{lemma}[Arithmetic regularity lemma]\label{arl}\cite{green}  Let $A \subset \F_2^n$ and $0 < \eps \leq 1$.  Then there exists a subspace $V$ of codimension $d = O_\eps(1)$ such that $A$ is $\eps$-regular on all but $\eps 2^d$ of the translates of $V$.
\end{lemma}

\begin{proof}  It will suffice to establish the claim with the weaker claim that $A$ is $O(\eps^{1/4})$-regular on all but $O(\sqrt{\eps} 2^d)$ of the translates of $V$, since one can simply shrink $\eps$ to obtain the original version of Lemma \ref{arl}.  

We apply Theorem \ref{sst} to the setting in Example \ref{fourier-struct}, with $f := 1_A$, and $F$ to be chosen later.  This gives us an integer $M = O_{F,\eps}(1)$ and a decomposition
\begin{equation}\label{adecom}
1_A = f_{\str} + f_\psd + f_\err
\end{equation}
where $f_\str$ is $(M,M)$-structured, $f_\psd$ is $1/F(M)$-pseudorandom, and $\|f_\err\|_H \leq \eps$.  The function $f_\str$ is a combination of at most $M$ characters, and thus there exists a subspace $V \subset \F_2^n$ of codimension $d \leq M$ such that $f_\str$ is constant on all translates of $V$.

We have
$$ \E_{x \in \F_2^n} |f_\err(x)|^2 \leq \eps = \eps 2^d |V|/|\F_2^n|.$$
Dividing $\F_2^n$ into $2^d$ translates $y+V$ of $V$, we thus conclude that we must have
\begin{equation}\label{xyv}
 \E_{x \in y+V} |f_\err(x)|^2 \leq \sqrt{\eps} 
\end{equation}
on all but at most $\sqrt{\eps} 2^d$ of the translates $y+V$.

Let $y+V$ be such that \eqref{xyv} holds, and let $\delta_{y+V}$ be the average of $A$ on $y+V$.  The function $f_\str$ equals a constant value on $y+V$, call it $c_{y+V}$.  Averaging \eqref{adecom} on $y+V$ we obtain
$$ \delta_{y+V} = c_{y+V} + \E_{x \in y+V} f_\psd(x) + \E_{x \in y+V} f_\err(x).$$
Since $f_\psd(x)$ is $1/F(M)$-pseudorandom, some simple Fourier analysis (expressing $1_{y+V}$ as an average of characters) shows that
$$\left|\E_{x \in y+V} f_\psd(x)\right| \leq \frac{2^n}{|V| F(M)} \leq \frac{2^M}{F(M)}$$
while from \eqref{xyv} and Cauchy-Schwarz we have
$$|\E_{x \in y+V} f_\err(x)| \leq \eps^{1/4}$$
and thus
$$ \delta_{y+V} = c_{y+V} + O\left( \frac{2^M}{F(M)} \right) + O( \eps^{1/4} ).$$
By \eqref{adecom} we therefore have
$$ 1_A(x) - \delta_{y+V} = f_\psd(x) + f_\err(x) + O\left( \frac{2^M}{F(M)} \right) + O( \eps^{1/4} ).$$
Now let $e_\xi$ be an arbitrary character.  By arguing as before we have
$$|\E_{x \in y+V} f_\psd(x) e_\xi(x)| \leq \frac{2^M}{F(M)}$$
and
$$|\E_{x \in y+V} f_\err(x) e_\xi(x)| \leq \eps^{1/4}$$
and thus
$$ \E_{x \in y+V} (1_A(x) - \delta_{y+V}) e_\xi(x) = O\left( \frac{2^M}{F(M)} \right) + O( \eps^{1/4} ).$$
If we now set $F(M) := \eps^{-1/4} 2^M$ we obtain the claim.
\end{proof}

For some applications of this lemma, see \cite{green}.  A decomposition in a similar spirit can also be found in \cite{bourgain-firstroth}, \cite{gk}.  The weak structure theorem for Reed-Muller codes was also employed in \cite{gt-finite}, \cite{green-montreal} (under the name of a \emph{Koopman-von Neumann type theorem}).

Now we obtain the Szemer\'edi regularity lemma itself.  Recall that if $G=(V,E)$ is a graph and $A, B$ are non-empty disjoint subsets of $V$, we say that the pair $(A,B)$ is \emph{$\eps$-regular} if for any $A' \subset A, B' \subset B$ with $|A'| \geq \eps |A|$ and $|B'| \geq \eps |B|$, the number of edges between $A'$ and $B'$ differs from $\delta_{A,B} |A'| |B'|$ by at most $\eps |A'| |B'|$, where $\delta_{A,B} = |E \cap (A \times B)|/|A| |B|$ is the edge density between $A$ and $B$.

\begin{lemma}[Szemer\'edi regularity lemma]\label{srl}\cite{szemeredi} Let $0 < \eps < 1$ and $m \geq 1$.  Then if $G=(V,E)$ is a graph with $|V|=n$ sufficiently large depending on $\eps$ and $m$, then there exists a partition $V = V_0 \cup V_1 \cup \ldots \cup V_{m'}$ with $m \leq m' \leq O_{\eps,m}(1)$ such that $|V_0| \leq \eps n$, $|V_1| = \ldots = |V_{m'}|$, and such that all but at most $\eps (m')^2$ of the pairs $(V_i, V_j)$ for $1 \leq i < j \leq m'$ are $\eps$-regular.
\end{lemma}

\begin{proof} It will suffice to establish the weaker claim that $|V_0| = O(\eps n)$, and all but at most $O(\sqrt{\eps} (m')^2)$ of the pairs $(V_i, V_j)$ are $O(\eps^{1/12})$-regular.  We can also assume without loss of generality that $\eps$ is small.

We apply Theorem \ref{sst} to the setting in Example \ref{prod} with $f := 1_E$ and $F$ to be chosen later.  This gives us an integer $M = O_{F,\eps}(1)$ and a decomposition
\begin{equation}\label{edecom}
1_E = f_{\str} + f_\psd + f_\err
\end{equation}
where $f_\str$ is $(M,M)$-structured, $f_\psd$ is $1/F(M)$-pseudorandom, and $\|f_\err\|_H \leq \eps$.  The function $f_\str$ is a combination of at most $M$ tensor products of indicator functions $1_{A_i \times B_i}$.  The sets $A_i$ and $B_i$ partition $V$ into at most $2^{2M}$ sets, which we shall refer to as \emph{atoms}.  If $|V|$ is sufficiently large depending on $M$, $m$ and $\eps$, we can then partition $V = V_0 \cup \ldots \cup V_{m'}$ with $m \leq m' \leq (m + 2^{2M}) / \eps$, $|V_0| = O(\eps n)$, $|V_1| = \ldots = |V_{m'}|$, and such that each $V_i$ for $1 \leq i \leq m'$ is entirely contained within an atom.  In particular $f_\str$ is constant on $V_i \times V_j$ for all $1 \leq i < j \leq m'$.  Since $\eps$ is small, we also have $|V_i| = \Theta( n/m')$ for $1 \leq i \leq m$.

We have
$$ \E_{(v,w) \in V \times V} |f_\err(v,w)|^2 \leq \eps$$
and in particular
$$ \E_{1 \leq i < j \leq m'} \E_{(v,w) \in V_i \times V_j} |f_\err(v,w)|^2 = O( \eps ).$$
Then we have
\begin{equation}\label{vivj}
 \E_{(v,w) \in V_i \times V_j} |f_\err(v,w)|^2 \leq \sqrt{\eps} 
\end{equation}
for all but $O(\sqrt{\eps} (m')^2)$ pairs $(i,j)$.

Let $(i,j)$ be such that \eqref{vivj} holds.   On $V_i \times V_j$, $f_\str$ is equal to a constant value $c_{ij}$.  Also, from the pseudorandomness of $f_\psd$ we have
\begin{align*}
|\sum_{(v,w) \in A' \times B'} f_\psd(v,w)| &\leq \frac{n^2}{F(M)} \\
&= O_{m,\eps,M}\left( \frac{|V_i| |V_j|}{F(M)} \right)
\end{align*}
for all $A' \subset V_i$ and $B' \subset V_j$.
By arguing very similarly to the proof of Lemma \ref{arl}, we can conclude that the edge density $\delta_{ij}$ of $E$ on $V_i \times V_j$ is
$$ \delta_{ij} = c_{ij} + O( \eps^{1/4} ) + O_{m,\eps,M}\left( \frac{1}{F(M)} \right)$$
and that
\begin{align*}
|\sum_{(v,w) \in A' \times B'} (1_E(v,w) - &\delta_{ij})| =  
\bigl(O( \eps^{1/4} ) \\
&+ O_{m,\eps,M}\left( \frac{1}{F(M)} \right)\bigr) |V_i| |V_j| 
\end{align*}
for all $A' \subset V_i$ and $B' \subset V_j$.  This implies that the pair $(V_i,V_j)$ is $O( \eps^{1/12} ) + O_{m,\eps,M}(1/F(M)^{1/3})$-regular.  The claim now follows by choosing $F$ to be a sufficiently rapidly growing function of $M$, which depends also on $m$ and $\eps$.  
\end{proof}

Similar methods can yield an alternate proof of the regularity lemma for hypergraphs \cite{gowers-hyper-4}, \cite{gowers-hyper}, \cite{rs}, \cite{rodl}; see \cite{tao-hyper}.  To oversimplify enormously, one works on higher product spaces such as $V \times V \times V$, and uses partial tensor products such as $(v_1,v_2,v_3) \mapsto 1_A(v_1) 1_E(v_2,v_3)$ as the structured objects.  The lower-order functions such as $1_E(v_2,v_3)$ which appear in the structured component are then decomposed again by another application of structure theorems (e.g. for $1_E(v_2,v_3)$, one would use the ordinary Szemer\'edi regularity lemma).  The ability to arbitrarily select the various functions $F$ appearing in these structure theorems becomes crucial in order to obtain a satisfactory hypergraph regularity lemma.  

See also \cite{afms} for another graph regularity lemma involving an arbitrary function $F$ which is very similar in spirit to Theorem \ref{sst}.  In the opposite direction, if one applies the weak structure theorem (Corollary \ref{wst}) to the product setting (Example \ref{prod}) one obtains a ``weak regularity lemma'' very close to that in \cite{fkannan}.

\Section{Structure and randomness in a measure space}

We have seen that the Hilbert space model for separating structure from randomness is satisfactory for many applications.  However, there are times when the ``$L^2$'' type of control given by this model is insufficient.  A typical example arises when one wants to decompose a function $f: X \to \R$ on a probability space $(X, \X, \mu)$ into structured and pseudorandom pieces, plus a small error.  Using the Hilbert space model (with $H = L^2(X)$), one can control the $L^2$ norm of (say) the structured component $f_\str$ by that of the original function $f$, indeed the construction in Theorem \ref{sst} ensures that $f_\str$ is an orthogonal projection of $f$ onto a subspace generated by some vectors in $S$.  However, in many applications one also wants to control the $L^\infty$ norm of the structured part by that of $f$, and if $f$ is non-negative one often also wishes $f_\str$ to be non-negative also.  More generally, one would like a \emph{comparison principle}: if $f, g$ are two functions such that $f$ dominates $g$ pointwise (i.e. $|g(x)| \leq f(x)$), and $f_\str$ and $g_\str$ are the corresponding structured components, we would like $f_\str$ to dominate $g_\str$.
One cannot deduce these facts purely from the knowledge that $f_\str$ is an orthogonal projection of $f$.  If however we have the stronger property that $f_\str$ is a \emph{conditional expectation} of $f$, then we can achieve the above objectives.  This turns out to be important when establishing structure theorems for \emph{sparse} objects, for which purely $L^2$ methods are inadequate; this was in particular a key point in the recent proof \cite{gt-primes} that the primes contained arbitrarily long arithmetic progressions.

In this section we fix the probability space $(X, \X, \mu)$, thus $\X$ is a $\sigma$-algebra on the set $X$, and $\mu: \X \to [0,1]$ is a probability measure, i.e. a countably additive non-negative measure.  In many applications one can assume that the $\sigma$-algebra $\X$ is finite, in which case it can be identified with a finite partition $X = A_1 \cup \ldots \cup A_k$ of $X$ into \emph{atoms} (so that $\X$ consists of all sets which can be expressed as the union of atoms).  

\begin{example}[Uniform distribution]\label{unif}  If $X$ is a finite set, $\X = 2^X$ is the power set of $X$, and $\mu(E) := |E|/|X|$ for all $E \subset X$ (i.e. $\mu$ is uniform probability measure on $X$), then $(X,\X,\mu)$ is a probability space, and the atoms are just singleton sets.
\end{example}

We recall the concepts of a \emph{factor} and of \emph{conditional expectation}, which will be fundamental to our analysis.

\begin{definition}[Factor]  A \emph{factor} of $(X,\X,\mu)$ is a triplet $\Y = (Y,\Y,\pi)$, where $Y$ is a set, $\Y$ is a $\sigma$-algebra, and $\pi: X \to Y$ is a measurable map.  If $\Y$ is a factor, we let $\B_\Y := \{ \pi^{-1}(E): E \in \Y \}$ be the sub-$\sigma$-algebra of $\X$ formed by pulling back $\Y$ by $\pi$.  A function $f: X \to \R$ is said to be \emph{$\Y$-measurable} if it is measurable with respect to $\B_\Y$.  If $f \in L^2(X,\X,\mu)$, we let $\E(f|Y) = \E(f|\B_Y)$ be the orthogonal projection of $f$ to the closed subspace $L^2(X,\B_Y,\mu)$ of $L^2(X,\X,\mu)$ consisting of $\Y$-measurable functions.  If $\Y = (Y,\Y,\pi)$ and $\Y' = (Y',\Y',\pi')$ are two factors, we let $\Y \vee \Y'$ denote the factor $(Y \times Y', \Y \otimes \Y', \pi \oplus \pi')$.
\end{definition}

\begin{example}[Colourings]  Let $X$ be a finite set, which we give the uniform distribution as in Example \ref{unif}.  Suppose we \emph{colour} this set using some finite \emph{palette} $Y$ by introducing a map $\pi: X \to Y$.  If we endow $Y$ with the discrete $\sigma$-algebra $\Y = 2^Y$, then $(Y,\Y,\pi)$ is a factor of $(X,\X,\mu)$.  The $\sigma$-algebra $\B_\Y$ is then generated by the \emph{colour classes} $\pi^{-1}(y)$ of the colouring $\pi$.  The expectation $\E(f|Y)$ of a function $f: X \to \R$ is then given by the formula $\E(f|Y)(x) := \E_{x' \in \pi^{-1}(\pi(x))} f(x')$ for all $x \in X$, where $\pi^{-1}(\pi(x))$ is the colour class that $x$ lies in.
\end{example}

In the previous section, the concept of structure was represented by a set $S$ of vectors.  In this section, we shall instead represent structure by a collection $\S$ of \emph{factors}.  We say that a factor $\Y$ has \emph{complexity} at most $M$ if it is the join $\Y = \Y_1 \vee \ldots \vee \Y_m$ of $m$ factors from $\S$ for some $0 \leq m \leq M$.  We also say that a function $f \in L^2(X)$ is \emph{$\eps$-pseudorandom} if we have $\|\E(f|\Y)\|_{L^2(X)} \leq \eps$ for all $\Y \in \S$.  We have an analogue of Lemma \ref{pseud-dec}:

\begin{lemma}[Lack of pseudorandomness implies energy increment]\label{pseud-dec-meas}  Let $(X,\X,\mu)$ and $\S$ be as above.  Let $f \in L^2(X)$ be such that $f-\E(f|\Y)$ is \emph{not} $\eps$-pseudorandom for some $0 < \eps \leq 1$ and some factor $\Y$.  Then there exists $\Y' \in \S$ such that
$\|\E(f|\Y \vee \Y') \|_{L^2(X)}^2 \geq \|\E(f|\Y)\|_{L^2(X)}^2 + \eps^2$.
\end{lemma}

\begin{proof} By hypothesis we have
$$ \| \E( f-\E(f|\Y) | \Y' ) \|_{L^2(X)}^2 \geq \eps^2$$
for some $\Y' \in \S$.  By Pythagoras' theorem, this implies that
$$ \| \E( f-\E(f|\Y) | \Y \vee \Y' ) \|_{L^2(X)}^2 \geq \eps^2.$$
By Pythagoras' theorem again, the left-hand side is $\| \E(f|\Y \vee \Y') \|_{L^2(X)}^2 - \| \E(f|\Y) \|_{L^2(X)}^2$, and the claim follows.
\end{proof}

We then obtain an analogue of Lemma \ref{wst2}:

\begin{lemma}[Weak structure theorem]\label{wst-mes} Let $(X,\X,\mu)$ and $\S$ be as above.  Let $f \in L^2(X)$ be such that $\|f\|_{L^2(X)} \leq 1$, let $\Y$ be a factor, and let $0 < \eps \leq 1$.  Then there exists a decomposition $f = f_\str + f_\psd$, where $f_\str = \E(f|\Y \vee \Y')$ for some factor $\Y'$ of complexity at most $1/\eps^2$, and $f_\psd$ is $\eps$-pseudorandom.
\end{lemma}

\begin{proof}  We construct factors $\Y_1, \Y_2, \ldots, \Y_m \in \S$ by the following algorithm:
\begin{itemize}
\item Step 0: Initialise $m=0$.
\item Step 1: Write $\Y' := \Y_1 \vee \ldots \vee \Y_m$, $f_\str:= \E(f|\Y \vee \Y')$, and $f_\psd := f - f_\str$.
\item Step 2: If $f_\psd$ is $\eps$-pseudorandom then \textbf{STOP}.  Otherwise, by Lemma \ref{pseud-dec-meas} we can find $\Y_{m+1}\in \S$ such that $\| \E(f|\Y \vee \Y' \vee \Y_{m+1}) \|_{L^2(X)}^2 \geq \| \E(f|\Y \vee \Y') \|_{L^2(X)}^2 + \eps^2$.
\item Step 3: Increment $m$ to $m+1$ and return to Step 1.
\end{itemize}
Since the ``energy'' $\|f_\str\|_{L^2(X)}^2$ ranges between $0$ and $1$ (by the hypothesis $\|f\|_{L^2(X)} \leq 1$) and increments by $\eps^2$ at each stage, we see that this algorithm terminates in at most $1/\eps^2$ steps.  The claim follows.
\end{proof}

Iterating this we obtain an analogue of Theorem \ref{sst}:

\begin{theorem}[Strong structure theorem]\label{sst2}  Let $(X,\X,\mu)$ and $\S$ be as above.  Let $f \in L^2(X)$ be such that $\|f\|_{L^2(X)} \leq 1$, let $\eps > 0$, and let $F: \Z^+ \to \R^+$ be an arbitrary function.  Then we can find an integer $M = O_{F,\eps}(1)$ and a decomposition \eqref{fdecomp} where $f_\str = \E(f|\Y)$ for some factor $\Y$ of complexity at most $M$, $f_\psd$ is $1/F(M)$-pseudorandom, and $f_\err$ has norm at most $\eps$.  
\end{theorem}

\begin{proof} Without loss of generality we may assume $F(M) \geq 2M$.  Also, it will suffice to allow $\Y$ to have complexity $O(M)$ rather than $M$.

We recursively define $M_0 := 1$ and $M_i := F(M_{i-1})^2$ for all $i \geq 1$.  We then recursively define factors $\Y_0, \Y_1, \Y_2, \ldots$ by setting $\Y_0$ to be the trivial factor, and then for each $i \geq 1$ using Lemma \ref{wst2} to find a factor $\Y'_i$ of complexity at most $M_i$ such that $f - \E(f|\Y_{i-1} \vee \Y'_i)$ is $1/F(M_{i-1})$-pseudorandom, and then setting $\Y_i := \Y_{i-1} \vee \Y'_i$.  By Pythagoras' theorem and the hypothesis $\|f\|_{L^2(X)} \leq 1$, the energy $\| \E(f|\Y_i) \|_{L^2(X)}^2$ is increasing in $i$, and is bounded between $0$ and $1$.  By the pigeonhole principle, we can thus find $1 \leq i \leq 1/\eps^2+1$ such that $\|\E(f|\Y_i)\|_{L^2(X)}^2 - \|\E(f|\Y_{i-1})\|_{L^2(X)}^2 \leq \eps^2$; by Pythagoras' theorem, this implies that $\|\E(f|\Y_i) - \E(f|\Y_{i-1})\|_{L^2(X)} \leq \eps$.  If we then set $f_\str := \E(f|\Y_{i-1})$, $f_\psd := f-\E(f|\Y_i)$, $f_\err := \E(f|\Y_i) - \E(f|\Y_{i-1})$, and $M := M_{i-1}$, we obtain the claim.
\end{proof}

This theorem can be used to give alternate proofs of Lemma \ref{arl} and Lemma \ref{srl}; we leave this as an exercise to the reader (but see \cite{tao-sz} for a proof of Lemma \ref{srl} essentially relying on Theorem \ref{sst2}).

As mentioned earlier, the key advantage of these types of structure theorems is that the structured component $f_\str$ is now obtained as a conditional expectation of the original function $f$ rather than merely an orthogonal projection, and so one has good ``$L^1$'' and ``$L^\infty$'' control on $f_\str$ rather than just $L^2$ control.  In particular, these structure theorems are good for controlling \emph{sparsely supported functions} $f$ (such as the normalised indicator function of a sparse set), by obtaining a densely supported function $f_\str$ which models the behaviour of $f$ in some key respects.  Let us give a simplified ``sparse structure theorem'' which is too restrictive for real applications, but which serves to illustrate the main concept.

\begin{theorem}[Sparse structure theorem, toy version]\label{sparse} Let $0 < \eps < 1$, let $F: \Z^+ \to \R^+$ be a function, and let $N$ be an integer parameter.  Let $(X,\X,\mu)$ and $\S$ be as above, and depending on $N$.  Let $\nu \in L^1(X)$ be a non-negative function (also depending on $N$) with the property that for every $M \geq 0$, we have the ``pseudorandomness'' property
\begin{equation}\label{ennu}
\| \E(\nu|\Y) \|_{L^\infty(X)} \leq 1 + o_M(1)
\end{equation}
for all factors $\Y$ of complexity at most $M$, where $o_M(1)$ is a quantity which goes to zero as $N$ goes to infinity for any fixed $M$.  Let $f: X \to \R$ (which also depends on $N$) obey the pointwise estimate $0 \leq f(x) \leq \nu(x)$ for all $x \in X$.  Then, if $N$ is sufficiently large depending on $F$ and $\eps$, we can find an integer $M = O_{F,\eps}(1)$ and a decomposition \eqref{fdecomp} where $f_\str = \E(f|\Y)$ for some factor $\Y$ of complexity at most $M$, $f_\psd$ is $1/F(M)$-pseudorandom, and $f_\err$ has norm at most $\eps$.  Furthermore, we have 
\begin{equation}\label{fstr}
0 \leq f_\str(x) \leq 1 + o_{F,\eps}(1)
\end{equation}
and
\begin{equation}\label{fstr-bong}
 \int_X f_\str\ d\mu = \int_X f\ d\mu.
\end{equation}
\end{theorem}

An example to keep in mind is where $X = \{1,\ldots,N\}$ with the uniform probability measure $\mu$, $\S$ consists of the $\sigma$-algebras generated by a single discrete interval $\{ n \in \Z: a \leq n \leq b \}$ for $1 \leq a \leq b \leq N$, and $\nu$ being the function $\nu(x) = \log N 1_A(x)$, where $A$ is a randomly chosen subset of $\{1,\ldots,N\}$ with $\P(x \in A) = \frac{1}{\log N}$ for all $1 \leq x \leq N$; one can then verify \eqref{ennu} with high probability using tools such as Chernoff's inequality.  Observe that $\nu$ is bounded in $L^1(X)$ uniformly in $N$, but is unbounded in $L^2(X)$. Very roughly speaking, the above theorem states that any dense subset $B$ of $A$ can be effectively ``modelled'' in some sense by a dense subset of $\{1,\ldots,N\}$, normalised by a factor of $\frac{1}{\log N}$; this can be seen by applying the above theorem to the function $f := \log N 1_B(x)$.

\begin{proof} We run the proof of Lemma \ref{wst-mes} and Theorem \ref{sst2} again.  Observe that we no longer have the bound $\|f\|_{L^2(X)} \leq 1$.  However, from \eqref{ennu} and the pointwise bound $0 \leq f \leq \nu$ we know that
\begin{align*}
 \| \E(f|\Y) \|_{L^2(X)} &\leq
\| \E(\nu|\Y) \|_{L^2(X)} \\
&\leq \| \E(\nu|\Y) \|_{L^\infty(X)} \\
&\leq 1 + o_M(1)
\end{align*}
for all $\Y$ of complexity at most $M$.  In particular, for $N$ large enough depending on $M$ we have
\begin{equation}\label{fbong}
\| \E(f|\Y) \|_{L^2(X)}^2 \leq 2
\end{equation}
(say).  This allows us to obtain an analogue of Lemma \ref{wst-mes} as before (with slightly worse constants), assuming that $N$ is sufficiently large depending on $\eps$, by repeating the proof more or less verbatim.  One can then repeat the proof of Theorem \ref{sst2}, again using \eqref{fbong}, to 
obtain the desired decomposition.  The claim \eqref{fstr} follows immediately from \eqref{ennu}, and \eqref{fstr-bong} follows since $\int_X \E(f|\Y)\ d\mu= \int_X f\ d\mu$ for any factor $\Y$.
\end{proof}

\begin{remark} In applications, one does not quite have the property \eqref{ennu}; instead, one can bound $\E(\nu|\Y)$ by $1 + o_M(1)$ outside of a small exceptional set, which has measure $o(1)$ with respect to $\mu$ and $\nu$.  In such cases it is still possible to obtain a structure theorem similar to Theorem \ref{sparse}; see \cite[Theorem 8.1]{gt-primes}, \cite[Theorem 3.9]{tao-gauss}, or \cite[Theorem 4.7]{ziegler}.  These structure theorems have played an indispensable role in establishing the existence of patterns (such as arithmetic progressions) inside sparse sets such as the prime numbers, by viewing them as dense subsets of sparse pseudorandom sets (such as the \emph{almost prime} numbers), and then appealing to a sparse structure theorem to model the original set by a much denser set, to which one can apply deep theorems (such as Szemer\'edi's theorem \cite{szemeredi}) to detect the desired pattern.  
\end{remark}

The reader may observe one slight difference between the concept of pseudorandomness discussed here, and the concept in the previous section.  Here, a function $f_\psd$ is considered pseudorandom if its conditional expectations $\E(f_{\psd}|\Y)$ are small for various structured $\Y$.  In the previous section, a function $f_\psd$ is considered pseudorandom if its correlations $\langle f_\psd, g \rangle_H$ were small for various structured $g$.  However, it is possible to relate the two notions of pseudorandomness by the simple device of using a structured function $g$ to generate a structured factor $\Y_g$.  In measure theory, this is usually done by taking the level sets $g^{-1}([a,b])$ of $g$ and seeing what $\sigma$-algebra they generate.  In many quantitative applications, though, it is too expensive to take \emph{all} of these the level sets, and so instead one only takes a finite number of these level sets to create the relevant factor.  The following lemma illustrates this construction:

\begin{lemma}[Correlation with a function implies non-trivial projection]  Let $(X,\X,\mu)$ be a probability space. Let $f \in L^1(X)$ and $g \in L^2(X)$ be such that $\|f\|_{L^1(X)} \leq 1$ and $\|g\|_{L^2(X)} \leq 1$.  Let $\eps > 0$ and $0 \leq \alpha < 1$, and let $\Y$ be the factor $\Y = (\R, \Y, g)$, where $\Y$ is the $\sigma$-algebra generated by the intervals $[(n+\alpha)\eps, (n+1+\alpha)\eps)$ for $n \in \Z$.  Then we have
$$ \| \E(f|\Y) \|_{L^2(X)} \geq |\langle f, g \rangle_{L^2(X)}| - \eps.$$
\end{lemma}

\begin{proof}  Observe that the atoms of $\B_\Y$ are generated by level sets $g^{-1}([(n+\alpha)\eps, (n+1+\alpha)\eps))$, and on these level sets $g$ fluctuates by at most $\eps$.  Thus
$$ \| g - \E(g|\Y) \|_{L^\infty(X)} \leq \eps.$$
Since $\|f\|_{L^1(X)} \leq 1$, we conclude
$$ \left|\langle f, g \rangle_{L^2(X)} - \langle f, \E(g|\Y) \rangle_{L^2(X)}\right| \leq \eps.$$
On the other hand, by Cauchy-Schwarz and the hypothesis $\|g\|_{L^2(X)} \leq 1$ we have
\begin{align*}
|\langle f, \E(g|\Y) \rangle_{L^2(X)}| &= |\langle \E(f|\Y), g \rangle_{L^2(X)}| \\
&\leq \| \E(f|\Y) \|_{L^2(X)}.
\end{align*}
The claim follows.
\end{proof}

This type of lemma is relied upon in the above-mentioned papers \cite{gt-primes}, \cite{tao-gauss}, \cite{ziegler} to convert pseudorandomness in the conditional expectation sense to pseudorandomness in the correlation sense.  In applications it is also convenient to randomise the shift parameter $\alpha$ in order to average away all boundary effects; see e.g. \cite[Lemma 3.6]{tao-norm}.

\Section{Structure and randomness via uniformity norms}

In the preceding sections, we specified the notion of structure (either via a set $S$ of vectors, or a collection $\S$ of factors), which then created a dual notion of pseudorandomness for which one had a structure theorem.  Such decompositions give excellent control on the structured component $f_\str$ of the function, but the control on the pseudorandom part $f_\psd$ can be rather weak.  There is an opposing approach, in which one first specifies the notion of pseudorandomness one would like to have for $f_\psd$, and then works as hard as one can to obtain a useful corresponding notion of structure.  In this approach, the pseudorandom component $f_\psd$ is easy to dispose of, but then all the difficulty gets shifted to getting an adequate control on the structured component.

A particularly useful family of notions of pseudorandomness arises from the \emph{Gowers uniformity norms} $\|f\|_{U^d(G)}$.  These norms can be defined on any finite additive group $G$, and for complex-valued functions $f: G \to \C$, but for simplicity let us restrict attention to a Hamming cube $G=\F_2^n$ and to real-valued functions $f: \F_2^n \to \R$.  (For more general groups and complex-valued functions, see \cite{tao-vu}.  For applications to graphs and hypergraphs, one can use the closely related \emph{Gowers box norms}; see \cite{gowers-hyper-4}, \cite{gowers-hyper}, \cite{lovasz-sz}, \cite{tao-gauss}, \cite{tao:montreal}, \cite{tao-vu}.)  In that case, the uniformity norm $\|f\|_{U^d(\F_2^n)}$ can be defined for $d \geq 1$ by the formula
$$ \|f\|_{U^d(\F_2^n)}^{2^d} := \E_{L: \F_2^d \to \F_2^n} \prod_{a \in \F_2^d} f(L(a))$$
where $L$ ranges over all affine-linear maps from $\F_2^d$ to $\F_2^n$ (not necessarily injective).  For instance, we have
\begin{align*}
\|f\|_{U^1(\F_2^n)} &= |\E_{x,h \in \F_2^n} f(x) f(x+h)|^{1/2} \\
&= |\E_{x \in \F_2^n} f(x)| \\
\|f\|_{U^2(\F_2^n)} &= |\E_{x,h,k \in \F_2^n} f(x) f(x+h) f(x+k) \\
&\quad \times f(x+h+k)|^{1/4} \\
&= |\E_{h \in \F_2^n} |\E_{x \in \F_2^n} f(x) f(x+h)|^2 |^{1/4}\\
\|f\|_{U^3(\F_2^n)} &= |\E_{x,h_1,h_2,h_3 \in \F_2^n} f(x) f(x+h_1) f(x+h_2)\\
&\quad \times f(x+h_3) f(x+h_1+h_2) f(x+h_1+h_3) \\
&\quad \times f(x+h_2+h_3) f(x+h_1+h_2+h_3)|^{1/8}.
\end{align*}
It is possible to show that the norms $\|\|_{U^d(\F_2^n)}$ are indeed a norm for $d \geq 2$, and a semi-norm for $d=1$; see e.g. \cite{tao-vu}.  These norms are also monotone in $d$:
\begin{equation}\label{mono}
 0 \leq \|f\|_{U^1(\F_2^n)} \leq \|f\|_{U^2(\F_2^n)} \leq \|f\|_{U^3(\F_2^n)} \leq \ldots \leq \|f\|_{L^\infty(\F_2^n)}.
 \end{equation}
The $d=2$ norm is related to the Fourier coefficients $\hat f(\xi)$ defined in \eqref{fourier-def} by the important (and easily verified) identity
\begin{equation}\label{fident}
\|f\|_{U^2(\F_2^n)} = (\sum_{\xi \in \F_2^n} |\hat f(\xi)|^4)^{1/4}.
\end{equation}
More generally, the uniformity norms $\|f\|_{U^d(\F_2^n)}$ for $d \geq 1$ are related to Reed-Muller codes of order $d-1$ (although this is partly conjectural for $d \geq 4$), but the relationship cannot be encapsulated in an identity as elegant as \eqref{fident} once $d \geq 3$.  We will return to this point shortly.

Let us informally call a function $f:\F_2^n \to \R$ \emph{pseudorandom of order $d-1$} if $\|f\|_{U^d(\F_2^n)}$ is small; thus for instance functions with small $U^2$ norm are \emph{linearly pseudorandom} (or \emph{Fourier-pseudorandom}, functions with small $U^3$ norm are \emph{quadratically pseudorandom}, and so forth.  It turns out that functions which are pseudorandom to a suitable order become negligible for the purpose of various multilinear correlations (and the higher the order of pseudorandomness, the more complex the multilinear correlations that become negligible).  This can be demonstrated by repeated application of the Cauchy-Schwarz inequality.  We give a simple instance of this:

\begin{lemma}[Generalised von Neumann theorem]  Let $T_1, T_2: \F^2_n \to \F_2^n$ be invertible linear transformations such that $T_1-T_2$ is also invertible. Then for any $f, g, h: \F^2_n \to [-1,1]$ we have
$$ |\E_{x,r \in \F_2^n} f(x) g(x+T_1 r) h(x+T_2 r)| \leq \|f\|_{U^2(\F_2^n)}.$$
\end{lemma}

\begin{proof} By changing variables $r' := T_2 r$ if necessary we may assume that $T_2$ is the identity map $I$.  We rewrite the left-hand side as
$$ |\E_{x \in \F_2^n} h(x) \E_{r \in \F_2^n} f(x-r) g(x+(T_1-I)r)|$$
and then use Cauchy-Schwarz to bound this from above by
$$ (\E_{x \in \F_2^n} |\E_{r \in \F_2^n} f(x-r) g(x+(T_1-I)r)|^2)^{1/2}$$
which one can rewrite as
$$ |\E_{x,r,r' \in \F_2^n} f(x-r) f(x-r') g(x+(T_1-I)r) g(x+(T_1-I)r')|^{1/2};$$
applying the change of variables $(y,s,h) := (x+(T_1-I)r,T_1 r,r-r')$, this can be rewritten as
$$ |\E_{y,h \in \F_2^n} g(y) g(y+(T_1-I)h) \E_{s \in \F_2^n} f(y+s) f(y+s+h)|^{1/2};$$
applying Cauchy-Schwarz, again, one can bound this by
$$ \left|\E_{y,h \in \F_2^n} |\E_{s \in \F_2^n} f(y+s) f(y+s+h)|^2\right|^{1/4}.$$
But this is equal to $\|f\|_{U^2(\F_2^n)}$, and the claim follows.
\end{proof}

For a more systematic study of such ``generalised von Neumann theorems'', including some weighted versions, see Appendices B and C of \cite{linear-primes}.

In view of these generalised von Neumann theorems, it is of interest to locate conditions which would force a Gowers uniformity norm $\|f\|_{U^d(\F_2^n)}$ to be small.  We first give a ``soft'' characterisation of this smallness, which at first glance seems too trivial to be of any use, but is in fact powerful enough to establish Szemer\'edi's theorem (see \cite{tao-ergodic}) as well as the Green-Tao theorem \cite{gt-primes}.  It relies on the obvious identity
$$ \|f\|_{U^d(\F_2^n)}^{2^d} = \langle f, {\mathcal D} f \rangle_{L^2(\F_2^n)}$$
where the \emph{dual function} ${\mathcal D} f$ of $f$ is defined as
\begin{equation}\label{dual}
 {\mathcal D} f(x) := \E_{L: \F_2^d \to \F_2^n; L(0)=x} \prod_{a \in \F_2^d \backslash \{0\}} f(L(a)).
\end{equation}

As a consequence, we have

\begin{lemma}[Dual characterisation of pseudorandomness]  Let $S$ denote the set of all dual functions ${\mathcal D} F$ with $\|F\|_{L^\infty(\F_2^n)} \leq 1$.  Then if $f: \F_2^n \to [-1,1]$ is such that
$\|f\|_{U^d(\F_2^n)} \geq \eps$ for some $0 < \eps \leq 1$, then we have $\langle f, g \rangle \geq \eps^{2^d}$ for some $g \in S$.
\end{lemma}

In the converse direction, one can use the \emph{Cauchy-Schwarz-Gowers inequality} (see e.g. \cite{gowers}, \cite{gt-primes}, \cite{linear-primes}, \cite{tao-vu}) to show that if $\langle f, g\rangle \geq \eps$ for some $g \in S$, then $\|f\|_{U^d(\F_2^n)} \geq \eps$.

The above lemma gives a ``soft'' way to detect pseudorandomness, but is somewhat unsatisfying due to the rather non-explicit description of the ``structured'' set $S$.  To investigate pseudorandomness further, observe that we have the recursive identity
\begin{equation}\label{recurse}
\|f\|_{U^d(\F_2^n)}^{2^d} = \E_{h \in \F_2^n} \| f f_h \|_{U^{d-1}(\F_2^n)}^{2^{d-1}}
\end{equation}
(which, incidentally, can be used to quickly deduce the monotonicity \eqref{mono}).  From this identity and induction we quickly deduce the modulation symmetry
\begin{equation}\label{fud}
\|fg\|_{U^d(\F_2^n)} = \|f\|_{U^d(\F_2^n)}
\end{equation}
whenever $g \in S_{d-1}(\F_2^n)$ is a Reed-Muller code of order at most $d-1$.  In particular, we see that $\|g\|_{U^d(\F_2^n)} = 1$ for such codes; thus a code of order $d-1$ or less is definitely \emph{not} pseudorandom of order $d$.  A bit more generally, by combining \eqref{fud} with \eqref{mono} we see that
$$ |\langle f, g \rangle_{L^2(\F_2^n)}| = \|fg\|_{U^1(\F_2^n)} \leq \|fg\|_{U^d(\F_2^n)} = \|f\|_{U^d(\F_2^n)}.$$
In particular, any function which has a large correlation with a Reed-Muller code $g \in S_{d-1}(\F_2^n)$ is not pseudorandom of order $d$.  It is conjectured that the converse is also true:

\begin{conjecture}[Gowers inverse conjecture for $\F_2^n$]  If $d \geq 1$ and $\eps > 0$ then there exists $\delta > 0$ with the following property: given any $n \geq 1$ and any $f: \F_2^n \to [-1,1]$ with $\|f\|_{U^d(\F_2^n)} \geq \eps$, there exists a Reed-Muller code $g \in S_{d-1}(\F_2^n)$ of order at most $d-1$ such that $|\langle f, g \rangle_{L^2(\F_2^n)}| \geq \delta$.
\end{conjecture}

This conjecture, if true, would allow one to apply the machinery of previous sections and then decompose a bounded function $f: \F_2^n \to [-1,1]$ (or a function dominated by a suitably pseudorandom function $\nu$) into a function $f_\str$ which was built out of a controlled number of Reed-Muller codes of order at most $d-1$, a function $f_\psd$ which was pseudorandom of order $d$, and a small error.  See for instance \cite{green-montreal} for further discussion.

The Gowers inverse conjecture is trivial to verify for $d=1$.  For $d=2$ the claim follows quickly from the identity \eqref{fident} and the Plancherel identity
$$ \|f\|_{L^2(\F_2^n)}^2 = \sum_{\xi \in \F_2^n} |\hat f(\xi)|^2.$$
The conjecture for $d=3$ was first established by Samorodnitsky \cite{sam}, using ideas from \cite{gowers-4} (see also \cite{gt:inverse-u3}, \cite{tao-vu} for related results).  The conjecture for $d > 3$ remains open; a key difficulty here is that there are a huge number of Reed-Muller codes (about $2^{\Omega(n^{d-1})}$ or so, compared to the dimension $2^n$ of $L^2(\F_2^n)$) and so we definitely do not have the type of orthogonality that one enjoys in the Fourier case $d=2$.  For related reasons, we do not expect any identity of the form \eqref{fident} for $d > 3$ which would allow the very few Reed-Muller codes which correlate with $f$ to dominate the enormous number of Reed-Muller codes which do not in the right-hand side.  

However, we can present some evidence for it here in the ``$99\%$-structured'' case when $\eps$ is very close to $1$.  Let us first handle the case when $\eps=1$:

\begin{proposition}[$100\%$-structured inverse theorem]  Suppose $d \geq 1$ and $f: \F_2^n \to [-1,1]$ is such that $\|f\|_{U^d(\F_2^n)}=1$.  Then $f$ is a Reed-Muller code of order at most $d-1$.
\end{proposition}

\begin{proof} We induct on $d$.  The case $d=1$ is obvious.  Now suppose that $d \geq 2$ and that the claim has already been proven for $d-1$.  If $\|f\|_{U^d(\F_2^n)}=1$, then from \eqref{recurse} we have
$$ \E_{h \in \F_2^n} \| f f_h \|_{U^{d-1}(\F_2^n)}^{2^{d-1}} = 1.$$
On the other hand, from \eqref{mono} we have $\| f f_h \|_{U^{d-1}(\F_2^n)} \leq 1$ for all $h$.  This forces $\| f f_h \|_{U^{d-1}(\F_2^n)} = 1$ for all $h$.  By induction hypothesis, $ff_h$ must therefore be a Reed-Muller code of order at most $d-2$ for all $h$.  Thus for every $h$ there exists a polynomial $P_h: \F_2^n \to \F_2$ of degree at most $d-2$ such that
$$ f(x+h) = f(x) (-1)^{P_h(x)}$$
for all $x, h \in \F_2^n$.  From this one can quickly establish by induction that for every $0 \leq m \leq n$, the function $f$ is a Reed-Muller code of degree at most $d-1$ on $\F_2^m$ (viewed as a subspace of $\F_2^n$), and the claim follows.
\end{proof}

To handle the case when $\eps$ is very close to $1$ is trickier (we can no longer afford an induction on dimension, as was done in the above proof).  We first need a rigidity result.

\begin{proposition}[Rigidity of Reed-Muller codes]\label{rigid}  For every $d \geq 1$ there exists $\eps > 0$ with the following property: if $n \geq 1$ and $f \in S_{d-1}(\F_2^n)$ is a Reed-Muller code of order at most $d-1$ such that $\E_{x \in \F_2^n} f(x) \geq 1-\eps$, then $f \equiv 1$.
\end{proposition}

\begin{proof} We again induct on $d$.  The case $d=1$ is obvious, so suppose $d \geq 2$ and that the claim has already been proven for $d-1$.  If $\E_{x \in \F_2^n} f(x) \geq 1-\eps$, then $\E_{x \in \F_2^n} |1-f(x)| \leq \eps$.  Using the crude bound $|1-ff_h| = O( |1-f| + |1-f_h| )$ we conclude that
$\E_{x \in \F_2^n} |1-ff_h(x)| \leq O(\eps)$, and thus
$$ \E_{x \in \F_2^n} ff_h(x) \geq 1-O(\eps)$$
for every $h \in \F_2^n$. But $ff_h$ is a Reed-Muller code of order $d-2$, thus by induction hypothesis we have $ff_h \equiv 1$ for all $h$ if $\eps$ is small enough.  This forces $f$ to be constant; but since $f$ takes values in $\{-1,+1\}$ and has average at least $1-\eps$, we have $f \equiv 1$ as desired for $\eps$ small enough.
\end{proof}

\begin{proposition}[$99\%$-structured inverse theorem]\cite{akklr}  For every $d \geq 1$ and $0 < \eps < 1$ there exists $0 < \delta < 1$ with the following property: if $n \geq 1$ and $f: \F_2^n \to [-1,1]$ is such that $\|f\|_{U^d(\F_2^n)} \geq 1-\delta$, then there exists a Reed-Muller code $g \in S_{d-1}(\F_2^n)$ such that $\langle f, g \rangle_{L^2(\F_2^n)} \geq 1-\eps$.
\end{proposition}

\begin{proof} We again induct on $d$.  The case $d=1$ is obvious, so suppose $d \geq 2$ and that the claim has already been proven for $d-1$.  Fix $\eps$, let $\delta$ be a small number (depending on $d$ and $\eps$) to be chosen later, and suppose $f: \F_2^n \to [-1,1]$ is such that $\|f\|_{U^d(\F_2^n)} \geq 1-\delta$.  We will use $o(1)$ to denote any quantity which goes to zero as $\delta \to 0$, thus $\|f\|_{U^d(\F_2^n)} \geq 1-o(1)$.  We shall say that a statement is true for \emph{most} $x \in \F_2^n$ if it is true for a proportion $1-o(1)$ of values $x \in \F_2^n$.

Applying \eqref{recurse} we have
$$ \E_{h \in \F_2^n} \| ff_h \|_{U^d(\F_2^n)} \geq 1-o(1)$$
while from \eqref{mono} we have $\| ff_h \|_{U^d(\F_2^n)} \leq 1$.  Thus we have $\|ff_h\|_{U^d(\F_2^n)} = 1-o(1)$ for all $h$ in a subset $H$ of $\F_2^n$ of density $1-o(1)$. Applying the inductive hypothesis, we conclude that for all $h \in H$ there exists a polynomial $P_h: \F_2^n \to \F_2$ of degree at most $d-2$ such that
$$ \E_{x \in \F_2^n} f(x) f(x+h) (-1)^{P_h(x)} \geq 1-o(1).$$
Since $f$ is bounded in magnitude by $1$, this implies for each $h \in H$ that
\begin{equation}\label{fxh}
f(x+h) = f(x) (-1)^{P_h(x)} + o(1)
\end{equation}
for most $x$.   For similar reasons it also implies that $|f(x)| = 1+o(1)$ for most $x$.

Now suppose that $h_1,h_2,h_3,h_4 \in H$ form an \emph{additive quadruple} in the sense that $h_1+h_2=h_3+h_4$.  Then from \eqref{fxh} we see that
\begin{equation}\label{fsh}
 f(x+h_1+h_2) = f(x) (-1)^{P_{h_1}(x)+P_{h_2}(x+h_1)} + o(1)
\end{equation}
for most $x$, and similarly
$$ f(x+h_3+h_4) = f(x) (-1)^{P_{h_3}(x)+P_{h_4}(x+h_3)} + o(1)$$
for most $x$.  Since $|f(x)| = 1+o(1)$ for most $x$, we conclude that
$$(-1)^{P_{h_1}(x)+P_{h_2}(x+h_1)-P_{h_3}(x)-P_{h_4}(x+h_3)} = 1+o(1)$$
for most $x$.  In particular, the average of the left-hand side in $x$ is $1-o(1)$.  Applying Lemma \ref{rigid} (and assuming $\delta$ small enough), we conclude that the left-hand side is \emph{identically} $1$, thus
\begin{equation}\label{pquad}
P_{h_1}(x)+P_{h_2}(x+h_1)=P_{h_3}(x)+P_{h_4}(x+h_3)
\end{equation}
for all additive quadruples $h_1+h_2=h_3+h_4$ in $H$ and all $x$.  

Now for any $k \in \F_2^n$, define the quantity $Q(k) \in \F_2$ by the formula
\begin{equation}\label{rk}
 Q(k) := P_{h_1}(0) + P_{h_2}(h_1)
\end{equation}
whenever $h_1,h_2 \in H$ are such that $h_1+h_2 \in H$.  Note that the existence of such an $h_1,h_2$ is guaranteed since most $h$ lie in $H$, and \eqref{pquad} ensures that the right-hand side of \eqref{rk} does not depend on the exact choice of $h_1,h_2$ and so $Q$ is well-defined.

Now let $x \in \F_2^n$ and $h \in H$.  Then, since most elements of $\F_2^n$ lie in $H$, we can find $r_1,r_2,s_1,s_2 \in H$ such that $r_1+r_2 = x$ and $s_1+s_2=x+h$.  From \eqref{fsh} we see that
$$ f(y+x) = f(y+r_1+r_2) = f(y) (-1)^{P_{r_1}(y) + P_{r_2}(y+r_1)} + o(1)$$
and
$$ f(y+x+h) = f(y+s_1+s_2) = f(y) (-1)^{P_{s_1}(y) + P_{s_2}(y+s_1)} + o(1)$$
for most $y$.  Also from \eqref{fxh}
$$ f(y+x+h) = f(y+x) (-1)^{P_h(y+x)} + o(1)$$
for most $y$.  Combining these (and the fact that $|f(y)| = 1+o(1)$ for most $y$) we see that
$$ (-1)^{P_{s_1}(y) + P_{s_2}(y+s_1) - P_{r_1}(y) - P_{r_2}(y+r_1) - P_h(y+x)} = 1+o(1)$$
for most $y$.  Taking expectations and applying Lemma \ref{rigid} as before, we conclude that
$$ P_{s_1}(y) + P_{s_2}(y+s_1) - P_{r_1}(y) - P_{r_2}(y+r_1) - P_h(y+x) = 0$$
for \emph{all} $y$.  Specialising to $y=0$ and applying \eqref{rk} we conclude that
\begin{equation}\label{phx}
P_h(x) = Q(x+h) - Q(x) = Q_h(x) - Q(x)
\end{equation}
for all $x \in \F_2^n$ and $h \in H$; thus we have succesfully ``integrated'' $P_h(x)$.  We can then extend $P_h(x)$ to all $h \in \F_2^n$ (not just $h \in H$) by viewing \eqref{phx} as a \emph{definition}.  Observe that if $h \in \F_2^n$, then $h=h_1+h_2$ for some $h_1,h_2 \in H$, and from \eqref{phx} we have
$$ P_h(x) = P_{h_1}(x) + P_{h_2}(x+h_1).$$
In particular, since the right-hand side is a polynomial of degree at most $d-2$, the left-hand side is also.  Thus we see that $Q_h-Q$ is a polynomial of degree at most $d-2$ for all $h$, which easily implies that $Q$ itself is a polynomial of degree at most $d-1$.  If we then set $g(x) := f(x) (-1)^{Q(x)}$, then from \eqref{fxh}, \eqref{phx} we see that for every $h \in H$ we have
$$ g(x+h) = g(x) + o(1)$$
for most $x$.  From Fubini's theorem, we thus conclude that there exists an $x$ such that $g(x+h) = g(x)+o(1)$ for most $h$, thus $g$ is almost constant.  Since $|g(x)| = 1+o(1)$ for most $x$, we thus conclude the existence of a sign $\epsilon \in \{-1,+1\}$ such that $g(x) = \epsilon + o(1)$ for most $x$.  We conclude that
$$ f(x) = \epsilon (-1)^{Q(x)} + o(1)$$
for most $x$, and the claim then follows (assuming $\delta$ is small enough).
\end{proof}

\begin{remark} The above argument requires $\|f\|_{U^d(\F_2^n)}$ to be very close to $1$ for two reasons.  Firstly, one wishes to exploit the rigidity property; and secondly, we implicitly used at many occasions the fact that if two properties each hold $1-o(1)$ of the time, then they jointly hold $1-o(1)$ of the time as well.  These two facts break down once we leave the ``$99\%$-structured'' world and instead work in a ``$1\%$-structured'' world in which various statements are only true for a proportion at least $\eps$ for some small $\eps$.  Nevertheless, the proof of the Gowers inverse conjecture for $d=2$ in \cite{sam} has some features in common with the above argument, giving one hope that the full conjecture could be settled by some extension of these methods.
\end{remark}

\begin{remark} The above result was essentially proven in \cite{akklr} (extending an argument in \cite{blm} for the linear case $d=2$), using a ``majority vote'' version of the dual function \eqref{dual}.
\end{remark}

\Section{Concluding remarks}

Despite the above results, we still do not have a systematic theory of structure and randomness which covers all possible applications (particularly for ``sparse'' objects).  For instance, there seem to be analogous structure theorems for random variables, in which one uses Shannon entropy instead of $L^2$-based energies in order to measure complexity; see \cite{tao-sz}.  In analogy with the ergodic theory literature (e.g. \cite{furst-book}), there may also be some advantage in pursuing \emph{relative} structure theorems, in which the notions of structure and randomness are all relative to some existing ``known structure'', such as a reference factor $\Y_0$ of a probability space $(X,\X,\mu)$.  Finally, in the iterative algorithms used above to prove the structure theorems, the additional structures used at each stage of the iteration were drawn from a fixed stock of structures ($S$ in the Hilbert space case, $\S$ in the measure space case).  In some applications it may be more effective to adopt a more \emph{adaptive} approach, in which the stock of structures one is using varies after each iteration.  A simple example of this approach is in \cite{tao-norm}, in which the structures used at each stage of the iteration are adapted to a certain spatial scale which decreases rapidly with the iteration.  I expect to see several more permutations and refinements of these sorts of structure theorems developed for future applications.

\Section{Acknowledgements}

The author is supported by a grant from the MacArthur Foundation, and by NSF grant CCF-0649473.  The author is also indebted to Ben Green for helpful comments and references.

\bibliographystyle{latex8}
\bibliography{latex8}

\end{document}